\documentclass{article}
\usepackage{amsfonts,amsmath,amssymb,cite}
\numberwithin{equation}{section}
\begin{document}
\title{An orthogonality relation for the Whittaker functions of the
second kind of imaginary order}
\author{Rados{\l}aw Szmytkowski\footnote{Email:
radek@mif.pg.gda.pl} \hspace{0em} and
Sebastian Bielski\footnote{Email: bolo@mif.pg.gda.pl} \\*[3ex]
Atomic Physics Division,
Department of Atomic Physics and Luminescence, \\
Faculty of Applied Physics and Mathematics,
Gda{\'n}sk University of Technology, \\
Narutowicza 11/12, PL 80--233 Gda{\'n}sk, Poland}
\date{\today}
\maketitle
\begin{abstract} 
An orthogonality relation for the Whittaker functions of the second
kind of imaginary order, $W_{\kappa,\mathrm{i}\mu}(x)$, with
$\mu\in\mathbb{R}$, is investigated. The integral
$\int_{0}^{\infty}\mathrm{d}x\:
x^{-2}W_{\kappa,\mathrm{i}\mu}(x)W_{\kappa,\mathrm{i}\mu'}(x)$ is
shown to be proportional to the sum
$\delta(\mu-\mu')+\delta(\mu+\mu')$, where $\delta(\mu\pm\mu')$ is
the Dirac delta distribution. The proportionality factor is found to
be $\pi^{2}/[\mu\sinh(2\pi\mu)
\Gamma(\frac{1}{2}-\kappa+\mathrm{i}\mu)
\Gamma(\frac{1}{2}-\kappa-\mathrm{i}\mu)]$. For $\kappa=0$ the
derived formula reduces to the orthogonality relation for the
Macdonald functions of imaginary order, discussed recently in the
literature. \\*[1ex]
\textbf{PACS:} 02.30.Gp \\*[1ex]
\textbf{MSC2010:} 33C15 \\*[1ex]
\textbf{Keywords:} Whittaker functions; orthogonal functions; 
Dirac delta distribution
\end{abstract}
%
%
\section{Introduction}
\label{I}
\setcounter{equation}{0}
The Whittaker functions $M_{\kappa,\lambda}(x)$ and
$W_{\kappa,\lambda}(x)$ \cite{Whit27,Buch69}, closely related to the
confluent hypergeometric functions, play an important role in various
branches of applied mathematics and theoretical physics, for instance
in fluid mechanics, scalar and electromagnetic diffraction theory or
atomic structure theory. This justifies the continuous effort in
studying properties of these functions and in gathering information
about them. It is the purpose of this brief paper to contribute to
the knowledge about the Whittaker function of the second kind.

The following double-integral formula:
\begin{eqnarray}
f(\mu) &=& \frac{\Gamma(\frac{1}{2}-\kappa+\mathrm{i}\mu)
\Gamma(\frac{1}{2}-\kappa-\mathrm{i}\mu)}{\pi^{2}}
\int_{0}^{\infty}\mathrm{d}x\:
\frac{W_{\kappa,\mathrm{i}\mu}(x)}{x^{2}}
\int_{0}^{\infty}\mathrm{d}\mu'\:\mu'\sinh(2\pi\mu')
W_{\kappa,\mathrm{i}\mu'}(x)f(\mu')
\nonumber \\
&& \hspace*{30em} (\mu>0),
\label{1.1}
\end{eqnarray}
valid under certain restrictions imposed on $f(\mu)$, was obtained by
Wimp \cite{Wimp64} as a particular case of a more general relation
involving the Meijer's $G$-function. From it one may infer that the
Whittaker functions of imaginary order, $W_{\kappa,\mathrm{i}\mu}(x)$
and $W_{\kappa,\mathrm{i}\mu'}(x)$, with $\mu,\mu'>0$, are orthogonal
on the positive real semi-axis with the weight $x^{-2}$ in the sense
of
\begin{equation}
\int\limits_{0}^{\infty}\mathrm{d}x\:
\frac{W_{\kappa,\mathrm{i}\mu}(x)
W_{\kappa,\mathrm{i}\mu'}(x)}{x^{2}}
=\frac{\pi^{2}}{\mu\sinh(2\pi\mu)
\Gamma(\frac{1}{2}-\kappa+\mathrm{i}\mu)
\Gamma(\frac{1}{2}-\kappa-\mathrm{i}\mu)}\,\delta(\mu-\mu')
\qquad (\mu,\mu'>0),
\label{1.2}
\end{equation}
where $\delta(\mu-\mu')$ is the Dirac delta function. Evidently, Eq.\
(\ref{1.2}) is to be understood in the distributional sense.

Of course, derivation of properties of the Whittaker's $W$-function
from those of the more general $G$-function is a perfectly valid
procedure. However, it is neither simple nor economical, as it
requires a good command of the rather complicated theory of the
Meijer's function. In this context, it seems natural to look for an
alternative, direct method of derivation of the relation (\ref{1.2}).
Such a method, making use of basic properties of the Whittaker
function and elementary facts from the theory of distributions only,
is presented in this paper. More precisely, in Sec.\ \ref{III} we
shall arrive at the relation
\begin{eqnarray}
\int\limits_{0}^{\infty}\mathrm{d}x\:
\frac{W_{\kappa,\mathrm{i}\mu}(x)
W_{\kappa,\mathrm{i}\mu'}(x)}{x^{2}}
&=& \frac{\pi^{2}}{\mu\sinh(2\pi\mu)
\Gamma(\frac{1}{2}-\kappa+\mathrm{i}\mu)
\Gamma(\frac{1}{2}-\kappa-\mathrm{i}\mu)}
[\delta(\mu-\mu')+\delta(\mu+\mu')]
\nonumber \\
&& \hspace*{20em} (\mu,\mu'\in\mathbb{R}),
\label{1.3}
\end{eqnarray}
which is slightly general than that in Eq.\ (\ref{1.2}) and reduces
to the latter for $\mu,\mu'>0$.

Throughout the rest of the work, unless otherwise stated, it is
assumed that $\kappa\in\mathbb{C}$, $\mu,\mu'\in\mathbb{R}$ and
$x\geqslant0$.
\section{Summary of relevant properties of the Whittaker functions of
the second kind of imaginary order} 
\label{II}
\setcounter{equation}{0}
Below we shall list these properties of the Whittaker functions of
the second kind of imaginary order, $W_{\kappa,\mathrm{i}\mu}(x)$,
which will prove to be helpful in Sec.\ \ref{III} for the derivation
of the relation (\ref{1.3}). The formulas presented below have been
extracted, with slight modifications whenever necessary, from the
invaluable collection by Magnus \emph{et al.\/} \cite{Magn66}.

The function $W_{\kappa,\mathrm{i}\mu}(x)$ satisfies the Whittaker
differential equation
\begin{equation}
\frac{\mathrm{d}^{2}F(x)}{\mathrm{d}x^{2}}
+\left(\frac{\mu^{2}+\frac{1}{4}}{x^{2}}+\frac{\kappa}{x}
-\frac{1}{4}\right)F(x)=0.
\label{2.1}
\end{equation}
The pair of the Whittaker functions of the first kind
\begin{equation}
M_{\kappa,\pm\mathrm{i}\mu}(x)
=x^{1/2\pm\mathrm{i}\mu}\mathrm{e}^{-x/2}
{}_{1}F_{1}({\textstyle\frac{1}{2}}-\kappa\pm\mathrm{i}\mu;
1\pm2\mathrm{i}\mu;x)
\label{2.2}
\end{equation}
also solves Eq.\ (\ref{2.1}). The functions
$M_{\kappa,\pm\mathrm{i}\mu}(x)$ and $W_{\kappa,\mathrm{i}\mu}(x)$
are not independent but are related through
\begin{equation}
W_{\kappa,\mathrm{i}\mu}(x)=\frac{\Gamma(2\mathrm{i}\mu)}
{\Gamma(\frac{1}{2}-\kappa+\mathrm{i}\mu)}M_{\kappa,-\mathrm{i}\mu}(x)
+\frac{\Gamma(-2\mathrm{i}\mu)}
{\Gamma(\frac{1}{2}-\kappa-\mathrm{i}\mu)}M_{\kappa,\mathrm{i}\mu}(x).
\label{2.3}
\end{equation}
Equation (\ref{2.3}) may serve as a definition of
$W_{\kappa,\mathrm{i}\mu}(x)$ in terms of 
$M_{\kappa,\pm\mathrm{i}\mu}(x)$.

For large positive values of $x$ the function
$W_{\kappa,\mathrm{i}\mu}(x)$ has the asymptotic representation
\begin{equation}
W_{\kappa,\mathrm{i}\mu}(x)\stackrel{x\to\infty}{\sim}
x^{\kappa}\mathrm{e}^{-x/2}
{}_{2}F_{0}({\textstyle\frac{1}{2}}-\kappa+\mathrm{i}\mu,
{\textstyle\frac{1}{2}}-\kappa-\mathrm{i}\mu;;-x^{-1}).
\label{2.4}
\end{equation}
From Eqs.\ (\ref{2.3}) and (\ref{2.2}) one finds that for small
positive values of $x$ the function $W_{\kappa,\mathrm{i}\mu}(x)$
behaves as
\begin{equation}
W_{\kappa,\mathrm{i}\mu}(x)\stackrel{x\to0+}{\sim}
x^{1/2}[A_{\kappa,\mathrm{i}\mu}\cos(-\mu\ln x)
+B_{\kappa,\mathrm{i}\mu}\sin(-\mu\ln x)][1+O(x^{-1})],
\label{2.5}
\end{equation}
where
\begin{equation}
A_{\kappa,\mathrm{i}\mu}=\frac{\Gamma(-2\mathrm{i}\mu)}
{\Gamma(\frac{1}{2}-\kappa-\mathrm{i}\mu)}
+\frac{\Gamma(2\mathrm{i}\mu)}
{\Gamma(\frac{1}{2}-\kappa+\mathrm{i}\mu)},
\label{2.6}
\end{equation}
\begin{equation}
\mathrm{i}B_{\kappa,\mathrm{i}\mu}
=\frac{\Gamma(-2\mathrm{i}\mu)}
{\Gamma(\frac{1}{2}-\kappa-\mathrm{i}\mu)}
-\frac{\Gamma(2\mathrm{i}\mu)}
{\Gamma(\frac{1}{2}-\kappa+\mathrm{i}\mu)}.
\label{2.7}
\end{equation}
\section{Orthogonality relation for the Whittaker functions of the
second kind of imaginary order} 
\label{III} 
\setcounter{equation}{0}
Consider two Whittaker functions $W_{\kappa,\mathrm{i}\mu}(x)$ and 
$W_{\kappa,\mathrm{i}\mu'}(x)$. According to what has been said in
Sec.\ \ref{II}, they obey the differential identities
\begin{equation}
\frac{\mathrm{d}^{2}W_{\kappa,\mathrm{i}\mu}(x)}{\mathrm{d}x^{2}}
+\left(\frac{\mu^{2}+\frac{1}{4}}{x^{2}}+\frac{\kappa}{x}
-\frac{1}{4}\right)W_{\kappa,\mathrm{i}\mu}(x)=0
\label{3.1}
\end{equation}
and
\begin{equation}
\frac{\mathrm{d}^{2}W_{\kappa,\mathrm{i}\mu'}(x)}{\mathrm{d}x^{2}}
+\left(\frac{\mu^{\prime\,2}+\frac{1}{4}}{x^{2}}+\frac{\kappa}{x}
-\frac{1}{4}\right)W_{\kappa,\mathrm{i}\mu'}(x)=0.
\label{3.2}
\end{equation}
We premultiply Eq.\ (\ref{3.1}) by $W_{\kappa,\mathrm{i}\mu'}(x)$,
Eq.\ (\ref{3.2}) by $W_{\kappa,\mathrm{i}\mu}(x)$, subtract and
integrate the result over $x$ from some $\xi>0$ to $\infty$. After
obvious movements, this gives
\begin{equation}
(\mu^{2}-\mu^{\prime\,2})\int\limits_{\xi}^{\infty}\mathrm{d}x\:
\frac{W_{\kappa,\mathrm{i}\mu}(x)W_{\kappa,\mathrm{i}\mu'}(x)}{x^{2}}
=\int\limits_{\xi}^{\infty}\mathrm{d}x\:
\left[W_{\kappa,\mathrm{i}\mu}(x)
\frac{\mathrm{d}^{2}W_{\kappa,\mathrm{i}\mu'}(x)}{\mathrm{d}x^{2}}
-W_{\kappa,\mathrm{i}\mu'}(x)
\frac{\mathrm{d}^{2}W_{\kappa,\mathrm{i}\mu}(x)}{\mathrm{d}x^{2}}
\right].
\label{3.3}
\end{equation}
Integrating the right-hand side by parts, we obtain
\begin{equation}
(\mu^{2}-\mu^{\prime\,2})\int\limits_{\xi}^{\infty}\mathrm{d}x\:
\frac{W_{\kappa,\mathrm{i}\mu}(x)W_{\kappa,\mathrm{i}\mu'}(x)}{x^{2}}
=\left[W_{\kappa,\mathrm{i}\mu}(x)
\frac{\mathrm{d}W_{\kappa,\mathrm{i}\mu'}(x)}{\mathrm{d}x}
-W_{\kappa,\mathrm{i}\mu'}(x)
\frac{\mathrm{d}W_{\kappa,\mathrm{i}\mu}(x)}{\mathrm{d}x}
\right]_{x=\xi}^{\infty}.
\label{3.4}
\end{equation}
Hence, after making use of the asymptotic property (\ref{2.4}), it
follows that
\begin{equation}
\int\limits_{0}^{\infty}\mathrm{d}x\:
\frac{W_{\kappa,\mathrm{i}\mu}(x)W_{\kappa,\mathrm{i}\mu'}(x)}{x^{2}}
=-\lim_{\xi\to0+}\frac{\displaystyle W_{\kappa,\mathrm{i}\mu}(\xi)
\frac{\mathrm{d}W_{\kappa,\mathrm{i}\mu'}(\xi)}{\mathrm{d}\xi}
-W_{\kappa,\mathrm{i}\mu'}(\xi)
\frac{\mathrm{d}W_{\kappa,\mathrm{i}\mu}(\xi)}{\mathrm{d}\xi}}
{\mu^{2}-\mu^{\prime\,2}}.
\label{3.5}
\end{equation}
To evaluate the limit on the right-hand side, we again exploit an
asymptotic property of the Whittaker function, this time the one in
Eq.\ (\ref{2.5}). This yields
\begin{eqnarray}
\int\limits_{0}^{\infty}\mathrm{d}x\:
\frac{W_{\kappa,\mathrm{i}\mu}(x)W_{\kappa,\mathrm{i}\mu'}(x)}{x^{2}}
&=& \lim_{\xi\to0+}
\left[\frac{\mu A_{\kappa,\mathrm{i}\mu}B_{\kappa,\mathrm{i}\mu'}
-\mu'B_{\kappa,\mathrm{i}\mu}A_{\kappa,\mathrm{i}\mu'}}
{\mu^{2}-\mu^{\prime\,2}}\sin(-\mu\ln\xi)\sin(-\mu'\ln\xi)\right.
\nonumber \\
&& +\,\frac{\mu A_{\kappa,\mathrm{i}\mu}A_{\kappa,\mathrm{i}\mu'}
+\mu'B_{\kappa,\mathrm{i}\mu}B_{\kappa,\mathrm{i}\mu'}}
{\mu^{2}-\mu^{\prime\,2}}\sin(-\mu\ln\xi)\cos(-\mu'\ln\xi)
\nonumber \\
&& -\,\frac{\mu B_{\kappa,\mathrm{i}\mu}B_{\kappa,\mathrm{i}\mu'}
+\mu'A_{\kappa,\mathrm{i}\mu}A_{\kappa,\mathrm{i}\mu'}}
{\mu^{2}-\mu^{\prime\,2}}\cos(-\mu\ln\xi)\sin(-\mu'\ln\xi)
\nonumber \\
&& -\left.\frac{\mu B_{\kappa,\mathrm{i}\mu}
A_{\kappa,\mathrm{i}\mu'}
-\mu'A_{\kappa,\mathrm{i}\mu}B_{\kappa,\mathrm{i}\mu'}}
{\mu^{2}-\mu^{\prime\,2}}\cos(-\mu\ln\xi)\cos(-\mu'\ln\xi)\right].
\label{3.6}
\end{eqnarray}
With a little trigonometry, this may be transformed into
\begin{eqnarray}
\int\limits_{0}^{\infty}\mathrm{d}x\:
\frac{W_{\kappa,\mathrm{i}\mu}(x)W_{\kappa,\mathrm{i}\mu'}(x)}{x^{2}}
&=& \frac{1}{2}\lim_{\xi\to0+}
\left\{(A_{\kappa,\mathrm{i}\mu}A_{\kappa,\mathrm{i}\mu'}
+B_{\kappa,\mathrm{i}\mu}B_{\kappa,\mathrm{i}\mu'})
\frac{\sin[-(\mu-\mu')\ln\xi]}{\mu-\mu'}\right.
\nonumber \\
&& +\,(A_{\kappa,\mathrm{i}\mu}A_{\kappa,\mathrm{i}\mu'}
-B_{\kappa,\mathrm{i}\mu}B_{\kappa,\mathrm{i}\mu'})
\frac{\sin[-(\mu+\mu')\ln\xi]}{\mu+\mu'}
\nonumber \\
&& +\,\frac{A_{\kappa,\mathrm{i}\mu}B_{\kappa,\mathrm{i}\mu'}
-B_{\kappa,\mathrm{i}\mu}A_{\kappa,\mathrm{i}\mu'}}{\mu-\mu'}
\cos[-(\mu-\mu')\ln\xi]
\nonumber \\
&& -\left.\frac{A_{\kappa,\mathrm{i}\mu}B_{\kappa,\mathrm{i}\mu'}
+B_{\kappa,\mathrm{i}\mu}A_{\kappa,\mathrm{i}\mu'}}{\mu+\mu'}
\cos[-(\mu+\mu')\ln\xi]\right\}.
\label{3.7}
\end{eqnarray}
It is known from the elementary theory of distributions (cf
\cite[Sec.\ 4.5]{Sned51}) that
\begin{equation}
\lim_{a\to\infty}\frac{\sin(ax)}{\pi x}
=\frac{1}{2\pi}\lim_{a\to\infty}\int_{-a}^{a}\mathrm{d}\eta\:
\mathrm{e}^{\mathrm{i}\eta x}=\delta(x)
\label{3.8}
\end{equation}
and also that (in the distributional sense)
\begin{equation}
\lim_{a\to\infty}\cos(ax)=0
\label{3.9}
\end{equation}
(the latter is the corollary from the Riemann--Lebesgue lemma). Since
for $\xi\to0+$ one has $-\ln\xi\to\infty$, application of the above
distributional relations and Eqs.\ (\ref{2.6}) and (\ref{2.7}) to
Eq.\ (\ref{3.7}) casts the latter into
{\small
\begin{eqnarray}
\int\limits_{0}^{\infty}\mathrm{d}x\:
\frac{W_{\kappa,\mathrm{i}\mu}(x)W_{\kappa,\mathrm{i}\mu'}(x)}{x^{2}}
&=& \pi\left[\frac{\Gamma(2\mathrm{i}\mu)\Gamma(-2\mathrm{i}\mu')}
{\Gamma(\frac{1}{2}-\kappa+\mathrm{i}\mu)
\Gamma(\frac{1}{2}-\kappa-\mathrm{i}\mu')}
+\frac{\Gamma(-2\mathrm{i}\mu)\Gamma(2\mathrm{i}\mu')}
{\Gamma(\frac{1}{2}-\kappa-\mathrm{i}\mu)
\Gamma(\frac{1}{2}-\kappa+\mathrm{i}\mu')}\right]\delta(\mu-\mu')
\nonumber \\
&& +\,\pi\left[\frac{\Gamma(2\mathrm{i}\mu)\Gamma(2\mathrm{i}\mu')}
{\Gamma(\frac{1}{2}-\kappa+\mathrm{i}\mu)
\Gamma(\frac{1}{2}-\kappa+\mathrm{i}\mu')}
+\frac{\Gamma(-2\mathrm{i}\mu)\Gamma(-2\mathrm{i}\mu')}
{\Gamma(\frac{1}{2}-\kappa-\mathrm{i}\mu)
\Gamma(\frac{1}{2}-\kappa-\mathrm{i}\mu')}\right]\delta(\mu+\mu')
\nonumber \\
&& \hspace*{20em} (\mu,\mu'\in\mathbb{R}).
\label{3.10}
\end{eqnarray}
}
This is the symmetric form of the orthogonality relation for the
Whittaker functions of the second kind of imaginary order. Upon
making use of the following known property of the Dirac delta
\cite[Sec.\ 4.4]{Sned51}:
\begin{equation}
g(\mu')\delta(\mu\mp\mu')=g(\pm\mu)\delta(\mu\mp\mu'),
\label{3.11}
\end{equation}
the relation (\ref{3.10}) may be rewritten compactly, although
unsymmetrically, as
\begin{equation}
\int\limits_{0}^{\infty}\mathrm{d}x\:
\frac{W_{\kappa,\mathrm{i}\mu}(x)W_{\kappa,\mathrm{i}\mu'}(x)}{x^{2}}
=\frac{2\pi|\Gamma(2\mathrm{i}\mu)|^{2}}
{\Gamma(\frac{1}{2}-\kappa+\mathrm{i}\mu)
\Gamma(\frac{1}{2}-\kappa-\mathrm{i}\mu)}
[\delta(\mu-\mu')+\delta(\mu+\mu')].
\label{3.12}
\end{equation}
This may be still simplified. Use of the formula
\begin{equation}
|\Gamma(2\mathrm{i}\mu)|=\sqrt{\frac{\pi}{2\mu\sinh(2\pi\mu)}}
\label{3.13}
\end{equation}
transforms Eq.\ (\ref{3.12}) into
\begin{eqnarray}
\int\limits_{0}^{\infty}\mathrm{d}x\:
\frac{W_{\kappa,\mathrm{i}\mu}(x)
W_{\kappa,\mathrm{i}\mu'}(x)}{x^{2}}
&=& \frac{\pi^{2}}{\mu\sinh(2\pi\mu)
\Gamma(\frac{1}{2}-\kappa+\mathrm{i}\mu)
\Gamma(\frac{1}{2}-\kappa-\mathrm{i}\mu)}
[\delta(\mu-\mu')+\delta(\mu+\mu')]
\nonumber \\
&& \hspace*{20em} (\mu,\mu'\in\mathbb{R}).
\label{3.14}
\end{eqnarray}
An alternative form of the above relation is obtained if one rewrites
the right-hand side with the use of the identity
\begin{equation}
\delta(\mu^{2}-\mu^{\prime\,2})=\frac{1}{2|\mu|}
[\delta(\mu-\mu')+\delta(\mu+\mu')].
\label{3.15}
\end{equation}
This yields
\begin{equation}
\int\limits_{0}^{\infty}\mathrm{d}x\:
\frac{W_{\kappa,\mathrm{i}\mu}(x)
W_{\kappa,\mathrm{i}\mu'}(x)}{x^{2}}
=\frac{2\pi^{2}}{\sinh(2\pi|\mu|)
\Gamma(\frac{1}{2}-\kappa+\mathrm{i}\mu)
\Gamma(\frac{1}{2}-\kappa-\mathrm{i}\mu)}\,
\delta(\mu^{2}-\mu^{\prime\,2})
\qquad (\mu,\mu'\in\mathbb{R}).
\label{3.16}
\end{equation}
If $\mu,\mu'>0$, then $\mu+\mu'\neq0$ and consequently in the
distributional sense it holds that
\begin{equation}
\delta(\mu+\mu')=0
\qquad (\mu,\mu'>0).
\label{3.17}
\end{equation}
In this case the relation (\ref{3.14}) reduces to
\begin{equation}
\int\limits_{0}^{\infty}\mathrm{d}x\:
\frac{W_{\kappa,\mathrm{i}\mu}(x)
W_{\kappa,\mathrm{i}\mu'}(x)}{x^{2}}
=\frac{\pi^{2}}{\mu\sinh(2\pi\mu)
\Gamma(\frac{1}{2}-\kappa+\mathrm{i}\mu)
\Gamma(\frac{1}{2}-\kappa-\mathrm{i}\mu)}\,\delta(\mu-\mu')
\qquad (\mu,\mu'>0).
\label{3.18}
\end{equation}

Concluding, we observe that if $\kappa$ is an integer or
half-integer, the product $\Gamma(\frac{1}{2}-\kappa+\mathrm{i}\mu)
\Gamma(\frac{1}{2}-\kappa-\mathrm{i}\mu)$ appearing in the relations
(\ref{3.14}), (\ref{3.16}) and (\ref{3.18}) may be expressed in terms
of elementary functions. Of particular interest is the case
$\kappa=0$ since then one has
\begin{equation}
W_{0,\mathrm{i}\mu}(x)=\sqrt{\frac{x}{\pi}}
K_{\mathrm{i}\mu}({\textstyle\frac{1}{2}}x),
\label{3.19}
\end{equation}
where $K_{\mathrm{i}\mu}(x)$ is the Macdonald function of imaginary
order used as a kernel in the Kontorovich--Lebedev transform.
Exploiting the identity
\begin{equation}
|\Gamma({\textstyle\frac{1}{2}+\mathrm{i}\mu})|
=\sqrt{\frac{\pi}{\cosh(\pi\mu)}},
\label{3.20}
\end{equation}
we find that the orthogonality relations (\ref{3.14}) and
(\ref{3.18}) go over into
\begin{equation}
\int\limits_{0}^{\infty}\mathrm{d}x\:
\frac{K_{\mathrm{i}\mu}(x)K_{\mathrm{i}\mu'}(x)}{x}
=\frac{\pi^{2}}{2\mu\sinh(\pi\mu)}
[\delta(\mu-\mu')+\delta(\mu+\mu')]
\qquad (\mu,\mu'\in\mathbb{R})
\label{3.21}
\end{equation}
and
\begin{equation}
\int\limits_{0}^{\infty}\mathrm{d}x\:
\frac{K_{\mathrm{i}\mu}(x)K_{\mathrm{i}\mu'}(x)}{x}
=\frac{\pi^{2}}{2\mu\sinh(\pi\mu)}\,\delta(\mu-\mu')
\qquad (\mu,\mu'>0),
\label{3.22}
\end{equation}
respectively. These relations have been derived by the present
authors in Ref.\ \cite{Szmy09} in the manner analogous to the one
which has led us above to Eqs.\ (\ref{3.14}) and (\ref{3.18}).
Somewhat earlier, Yakubovich \cite{Yaku06} and Passian \emph{et
al.\/} \cite{Pass09} proved the validity of the relation (\ref{3.22})
exploiting more involved mathematical techniques.
%
%

%
\end{document}